\documentclass[12pt,a4,amsmath, amsthm, amssymb]{amsart}

\newtheorem{Lemma}{LEMMA}
\newtheorem{Theorem}[Lemma]{Theorem}
\newtheorem{Proposition}[Lemma]{Proposition}
\newtheorem{Corollary}[Lemma]{Corollary}

\newtheorem{remark}[Lemma]{Remark}
\newtheorem{definition}[Lemma]{Definition}
\newtheorem{example}[Lemma]{Example}

\newtheorem{Fact}[Lemma]{Fact}

\def\bt{\begin{Theorem}}
\def\et{\end{Theorem}}
\def\bl{\begin{Lemma}}
\def\el{\end{Lemma}}
\def\bp{\begin{Proposition}}
\def\ep{\end{Proposition}}
\def\bcor{\begin{Corollary}}
\def\ecor{\end{Corollary}}
\def\bpf{\begin{proof}}
\def\epf{\end{proof}}

\def\brem{\begin{remark}}
\def\erem{\end{remark}}

\def\bedef{\begin{definition}}
\def\endef{\end{definition}}

\def\beg{\begin{example}}
\def\eeg{\end{example}}

\def\bef{\begin{Fact}}
\def\eef{\end{Fact}}

\def\bc{\begin{center}}
\def\ec{\end{center}}
\def\noi{\noindent}

\def\vsq{\vskip .25cm}
\def\beq{\begin{equation}}
\def\eeq{\end{equation}}
\def\beqarray{\begin{eqnarray*}}
\def\eeqarray{\end{eqnarray*}}
\def\<{\leftangle}
\def\>{\rightangle}
\def\({\left(}
\def\){\right)}
\def\f{\varphi}

\def\<{\langle}
\def\>{\rangle}
\def\q{\quad}

\def\a{\alpha}
\def\b{\beta}

\def\p{\pi}
\def\d{\delta}
\def\h{\hbox}

\def\t{\tau}

\def\l{\lambda}
\def\e{\varepsilon}
\def\s{\sigma}
\def\O{\Omega}

\def\w.r.t.{with respect to}
\def\R{{\mathbb{R}}}
\def\N{{\mathbb{N}}}

\def\C{{\mathbb{C}}}

\def\K{{\mathbb{K}}}

\def\bq{\begin{quote}}
\def\eq{\end{quote}}

\def\bit{\begin{itemize}}
\def\eit{\end{itemize}}
\def\i{\item}
\def\ben{\begin{enumerate}}
\def\een{\end{enumerate}}

\def\ds{\displaystyle}

\begin{document}

\title[ A Discrete Regularization  Method  ]
{A Discrete Regularization Method for  \\  Ill-Posed Operaror Equations}
\author{M.T. Nair }
\address{Department of Mathematics, I.I.T. Madras,
Chennai-600 036, INDIA}
\email{mtnair@iitm.ac.in}
\today
\maketitle


\begin{abstract}
Discrete regularization methods are often applied for obtaining stable approximate solutions for ill-posed operator equations $Tx=y$, where $T: X\to Y$ is a bounded operator between Hilbert spaces with non-closed range $R(T)$ and $y\in R(T)$.  Most of the existing such methods  involve finite rank bounded projection operators on either the domain space $X$ or on codomain space $Y$ or on both. In this paper, we propose a discrete regularization based on finite rank projection-like operators  on some subspace of the codomain space such that their  ranges need not be subspaces of the codomain space. This method not only incudes some of the exisiting projection based methods but also a quadrature based collocation method considered by the author in \cite{mtn-acm} for integral equations of the firt kind.

%
\end{abstract}


\section{Introduction and Preparatory Results}

Consider the problem of solving the operator equation 
\beq\label{eq-1} Tx=y,\eeq
where  $T: X\to Y$ is a bounded linear operator between Hilbert spaces  $X$ and $Y$,  and $y\in Y$. 
We may observe that if $R(T)$, the range of $T$, is not closed in $Y$, then  $T$ cannot have a continuous inverse, and hence, the probelm of solving the  equation (\ref{eq-1}) is {\it ill-posed},  in the sense that  (\ref{eq-1}) need not have a solution for a given $y\in Y$, and even if it has a unique soolution,  the solution, it does not depend continuously on the data $y$.
In such case,  to obtain a  stable approximate solution, some regularization method has to be used. One such method, in the setting of infinite dimesional setting, is {\it Tikhonov regularization} (see \cite{tikh-ars, Gr-Tikh, ehn, mtn-linop}).   

In this paper we consider  a  discrete regularization method based on a finite rank projection-like operator $\pi_n$ when $y\in R(T)$. 
To prove the convergence of the method and to obtain the error estimate,  we shall make use of the error in Tikhonov regularization  and some assumptions on the sequence of operators which converge to $T^*T$. The special cases of the suggested method include various projection-like methods and also  a quadrature based collocation method consdered in  \cite{mtn-acm, mtn-sp-jc}.

Since the error estimate and the convergence results heavily depend on the Tikhonov regularization, let us recall some results pertaining to it. Recall that 
 the regularized approximation is defined as the solution of the well-posed equation 
\beq\label{eq-Tikh}(T^*T+\a I) x_\a = T^*y\eeq
for $\a>0$.  It is known that, if $y\in R(T)+R(T)^\perp$, then there exists a unique $x^\dagger \in X$ such that 
$$\|x^\dagger\| = \inf\{\|x\|: \|Tx-y\|\leq \|Tu-y\|\,\forall\, u\in X\},$$
and it is the unique element in $N(T)^\perp$ such that 
\beq\label{eq-normal} T^*Tx^\dagger =T^*y.\eeq
It is also known that the map $T^\dagger: R(T)+R(T)^\perp\to X$, called the {\it genrealized inverse} of $T$, which maps $y$ to  $x^\dagger:= T^\dagger y$ is a closed operator, and it is continuous if and only if $R(T)$ is closed.

Suppose $y\in R(T)+R(T)^\perp$  and $x^\dagger =T^\dagger y$.  
Then from (\ref{eq-Tikh}) and (\ref{eq-normal}), we obtain  
$$x^\dagger  - x_\a = \a(T^*T+\a I)^{-1}x^\dagger.$$
Since $x^\dagger \in N(T)^\perp$, and since $R(T^*T)$ is a dense subspace of $N(T)^\perp$, it can be shown by using a stadard result in Functional Anaysis (cf. \cite{mtn-fa}) that  
$$\|x^\dagger  - x_\a\|\to 0\q\h{ as}\q \a\to 0.$$
However, for obtianing an estimate for $\|x^\dagger  - x_\a\|$,  additional assumptions has to be imposed on the $x^\dagger$. In this regard, we have the following result  (see \cite{ mtn-linop} for its proof).

\bt\label{Th-Tikh}   
Suppose $y\in R(T)+R(T)^\perp$  and $x^\dagger =T^\dagger y$. Let  $x_\a$ be as in  (\ref{eq-Tikh}).
 If $x^\dagger = \f(T^*T)u$\,  for some $u\in X$, where  $\f: (0,\infty)\to (0, \infty)$ is a continuous function such that 
$\f(\a)\to 0$ as $\a\to 0$ and 
\beq\label{source} \sup_{\l>0}\frac{\a\f(\l)}{\l +\a}\leq c_0\f(\a)\eeq
for some  $c_0>0$, then
$$\|x^\dagger  - x_\a\|\leq c_0\|u\|\f(\a).$$
\et
Particular choices of $\f$ in (\ref{source})  are  
\ben
\i[(i)]  $\f(\l):=\l^\nu$ for some $\nu\in (0, 1]$ or 
\i[(ii)]  $\f(\l):= [\log(1/\l)]^{-p}$ for some $p>0.$
\een
The   cases (i) and (ii) above are useful when the problem (\ref {eq-1}) is mildly ill-posed and severely ill-posed, respectively, (see \cite{mair, mtn-linop})  (c.f. Hohag\cite{hohag}). 

When $T$ is a general bounded linear operator, the oprator  $\f(T^*T)$   in  Theorem \ref{Th-Tikh},  is  defined via the {\it spectal theorem}  for the self-adjoint (positive) operator $T^*T$ (cf. \cite{mtn-fa})hat is, 
$$\f(T^*T) = \int_0^a\l dE_\l,\q a>\|T\|^2,$$
where $\{E_\l: \l \in [0, a]\}$ is the {\it spectral resolution of identity} corresponding to the self adjoint operator $T^*T$ (see \cite{yoshida, mtn-fa}).
If $T$ is a compact oprator of infite rank, then one may use the  singular value decomposition (see \cite{mtn-fa}) of $T$, that is, 
\beq\label{svd} 
Tx = \sum_{n=1}^\infty \s_n\<x, u_n\> v_n,\q x\in X,
\eeq
to define $\f(T^*T)$, that is, 
$$\f(T^*T)x= \sum_{n=1}^\infty \f(\s_n)\<x, u_n\> v_n,\q x\in X.$$
Here, $\{u_n: n\in \N\}$ and $\{v_n: n\in \N\}$ are orthonormal bases of $N(T)^\perp$ and $N(T^*)^\perp$, respectively and the {\it singular values}  $\s_n$'s of $T$ are positive real numbers such that $\s_n\to 0$ as $n\to \infty$. In this case, we have 
\beq\label{gen-sol} 
T^\dagger y = \sum_{n=1}^\infty \frac{\<y, v_n\>}{\s_n} u_n,\q y\in R(T)+R(T)^\perp,
\eeq
and 
 the assumption  $x^\dagger = \f(T^*T)u$ is equivalent to the requirement   
$$\sum_{n=1}^\infty \frac{|\<x^\dagger, u_n\>|^2}{\f(\s_n)^2}<\infty. $$
(see (cf. \cite{mtn-linop}). It is also to be obseved that  if $R(T)$ is of finite rank, say ${\rm rank}(T)=r$, then  (\ref{svd}) and ({\ref{gen-sol}) take the forms 
$$Tx = \sum_{n=1}^r \s_n\<x, u_n\> v_n,\q x\in X,$$
and 
$$T^\dagger y = \sum_{n=1}^r \frac{\<y, v_n\>}{\s_n} u_n,\q y\in R(T)+R(T)^\perp,$$
respectively, and in that case,  
$$\|T^\dagger\| = \frac{1}{\ds \min_{1\leq n\leq r} \s_n}.$$

Suppose the data is noisy, say  $\tilde y$ in place of $y$ with $\|y-\tilde y\|\leq \d$. If  $\tilde x_\a$ is the corresponding Tikhonov regularized solution, that is, 
\beq\label{eq-Tikh-noisy} 
(T^*T+\a I) \tilde x_\a = T^*\tilde y
\eeq
 and if $R(T)$ is not closed,  then the family $(\|\tilde x_\a\|)$ need not be even bounded (see \cite{ehn, mtn-linop}).
In such case, 
one has to choose the {\it regularization parameter} $\a$ depending on $(\d, \tilde y)$ in such a way that 
$$\|x^\dagger - \tilde x_\a\|\to 0\q\h{as}\q \d\to 0.$$
There are many parameter choice available in the literature (see \cite{ehn, mtn-linop}) and some of the references there in).

In Section 2 we introduce the proposed discrete regularization method and derive a general error estimate which motivates the type of assumptions to be made for establishing the convergence. 
In Section 3 we consider some of the special cases of the general setting of Section 2, and in Section 4 we consider the case when the data is noisy. In Section 5 further regularization of the discretised  equation  is considered and in the final section, Section 6, we discuss the method in the context of integral equations of the first kind,  inclduing the one considered in  \cite{mtn-acm, mtn-sp-jc}. 

\section{The Method and a general Error Estimate} 

Throughout the paper $T: X\to Y$ is a bounded linear operaror between Hilbert spaces $X$ and $Y$ over the space scalar filed $\K$ which is either $\R$ or $\C$. We also assume the following:
\ben
\i[(a)]   $Z$ is a  subspace of $Y$ such that $R(T)\subseteq Z$ and $y\in R(T)$.
\i[(b)] For each $n\in \N$, $Y_n$ is a finite dimensional inner product space over $\K$.   
\i[(c)] For each $n\in \N$,   $\pi_n: Z\to Y_n$ is a  linear operator such that $\p_nT: X\to  Y_n$ is a bounded linear operator.  
\een

It is to be observed that, though $\pi_n$ is a finite reank operator, we do not assume it  to be a bounded linear operaror.  

Since $y\in R(T)$, we have  
Then we also have   
\beq\label{eq-2} Tx^\dagger = y \q\h{and}\q T_nx^\dagger = y_n,\eeq
where   $x^\dagger = T^\dagger y$.
However, since $T^\dagger$ is an unbounded operator, $x^\dagger$ cannot be recovered in a stable manner.
Since $\pi_nT$ is of finite rank, a natural choice would be to consider the minimum norm solution $x_n^\dagger$ of 
\beq\label{eq-3}  \pi_nTx = \pi_ny.
\eeq
We may recall that $x_n^\dagger$ is the unique element in $N(\pi_nT)^\perp$ such that 
$$\|x_n^\dagger\| = \inf\{\|u\|: u\in X,\, \pi_nTx=\pi_ny\}.$$
We denote 
$$T_n:= \pi_nT\q\h{and}\q y_n=\pi_ny.$$
Note that, since  $T_n: X\to Y_n$ is of finite rank,  $N(T_n)^\perp=R(T_n^*)$ and hence, $N(T_n)^\perp$ is finite dimensional so that 
$x_n^\dagger$ is the unique solution of the equation
$${T_n}_{|_{N(T_n)^\perp}} x=y_n,$$
and it can be also be represented as 
$$x_n^\dagger = T_n^*v_n,$$
where $v_n\in Y_n$ is any solution of   
\beq\label{eq-comp-1}  
T_nT_n^* v_n =  y_n.
\eeq
Recall that $T_nT_n^*$ is an  operator from $Y_n$ into itself. Therefore, the equation (\ref{eq-comp-1}) corresponds to solving a matrix equation.


\bt\label{Th-1}
Let $\e_n>0$ be such that 
$\|T^*T- T_n^*T_n\|\leq \e_n$. Then for every $\a>0$, 
$$\|x^\dagger -x_n^\dagger \| \leq  \Big(1+\frac{\e_n}{\a}\Big) \|x^\dagger -x_\a\|.$$
In particular,  the following results hold.
\ben
\i  \, $\|x^\dagger -x_n^\dagger \| \leq 2  \|x^\dagger -x_{\e_n}\|.$
\i \, If $\|T^*T- T_n^*T_n\|\to 0$ as $n\to\infty$, then $\|x^\dagger -x_n^\dagger \|\to 0$ as $n\to\infty$.
\een
\et

\bpf
Since $\pi_n y\in R(T_n)$, we have 
\beq\label{eq-3} T_nx_n^\dagger = \pi_ny.\eeq
From  (\ref{eq-2}) and (\ref{eq-3}), we  obtain 
$$x^\dagger -x_n^\dagger \in N(T_n).$$
Since $x_n^\dagger\in N(T_n)^\perp$, it follows that 
$$x_n^\dagger=Q_nx^\dagger,$$
where $Q_n: X\to X$ is the orthogonal projection onto $N(T_n)^\perp$. Thus, 
$$\|x^\dagger -x_n^\dagger \|\leq \|x^\dagger - u\|\q\forall\, u\in N(T_n)^\perp.$$
Since $R(T_n)$ is closed, we have   
$$ N(T_n)^\perp = N(T_n^*T_n)^\perp=R(T_n^*T_n).$$ 
Also, 
 $x_\a = T^*Tz_\a,$ where 
$$z_\a = (T^*T+\a I)^{-1}x^\dagger = \frac{1}{\a}(x^\dagger-x_\a).$$
Thus, 
\beqarray
\|x^\dagger -x_n^\dagger \| &\leq & \|x^\dagger -x_\a\| + \|x_\a - T_n^*T_nz_\a\|\\
 & =  & \|x^\dagger -x_\a\| + \|(T^*T-  T_n^*T_n)z_\a\|\\
 & =  & \|x^\dagger -x_\a\| + \frac{1}{\a}\|(T^*T- T_n^*T_n) (x^\dagger-x_\a) \|\\
& \leq   & (1+\frac{\e_n}{\a}) \|x^\dagger -x_\a\|,
\eeqarray
where $\e_n>0$ such that 
$\|T^*T- T_n^*T_n\|\leq \e_n.$  The particular cases are obvious. 
\epf


\section{Special cases}

By Theorem \ref{Th-1},  the convergence of the method is guaranteed if 
\beq\label{conv-op-1} 
\|T^*T- T_n^*T_n\|\to 0\q\h{as}\q n\to\infty.
\eeq
We shall specify specify some  conditions under which we have this convergence. We shall also  make use of the following result  from functional analysis (cf.  \cite{mtn-fa}).

\bp\label{Prop-fa}
Let $X_1, X_2,  X_3$ be normed linear sapces.  If $X_2$ is a Banach space, $A: X_1\to X_2$ is a compact operator,  $B: X_2\to X_3$ is a bounded linear operaror and $(B_n)$ is a sequence of bounded linear operators from $X_2$ to $X_3$  such that $\|B_nx-Bx\|\to 0$ as $n\to\infty$ for each $x\in X_2$, then 
$$\|(B_n-B)A\|\to 0\q\h{as}\q n\to\infty.$$
\ep

\bt\label{Th-special-1} 
Suppose  that for each $n\in \N$, $Y_n$ is a  subspace of $Y$.   Then 
$$\|T^*T-T_n^*T_n\| \leq (\|T\| +  \|T_n\|) \|(I-\pi_n)T\|.$$
In particular, if  
\beq\label{conv-op-2}\|(I-\pi_n)T\|\to 0\q\h{as}\q n\to\infty,
\eeq
then  $(\|T_n\|)$ is bounded and 
$\|T^*T- T_n^*T_n\|\to 0$ {as} $n\to\infty.$
\et

\bpf
Let $x\in X$. Then we have 
\beqarray
\<(T^*T-T_n^*T_n)x, x\> &=& \<Tx, Tx\> - \<T_nx, T_n\>\\
&=& \<(T-T_n)x, Tx\> + \<T_nx, (T-T_n)x\>\\
&=& \<(I-\pi_n)Tx, Tx\> + \<\pi_nTx, (I-\pi_n)Tx\>.
\eeqarray
Hence,
$$|\<(T^*T-T_n^*T_n)x, x\>|\leq   (\|T\| +  \|T_n\|) \|(I-\pi_n)T\| \|x\|^2.$$
Since $T^*T-T_n^*T_n$ is a  self adjoint operator, we know that  (cf. \cite{mtn-fa}) 
$$\|T^*T-T_n^*T_n\| = \sup_{\|x\|=1} |\<(T^*T-T_n^*T_n)x, x\>|.$$
 Thus,
$$\|T^*T-T_n^*T_n\| \leq   (\|T\| +  \|T_n\|)  \|(I-\pi_n)T\|.$$

Now, assume the condition (\ref{conv-op-2}) . Then   $(\|\pi_nT\|)$ is bounded, say $\|\pi_nT\|\leq c$ for all $n\in \N$. Hence, we have 
$$\|T^*T-T_n^*T_n\| \leq  (\|T\| + c) \|(I-\pi_n)T\|\to 0\q\h{as}\q n\to\infty.$$
\epf

\brem{\rm 
We may also observe that if $\pi_n$ are orthogonal projections on $Y$ with $R(\pi_n)=Y_n$, then 
\beqarray
\|T^*T- T_n^*T_n\| &=& \|T^*T - T^*\pi_nT\| \\ 
&=& \|T^*(I-\pi_n)T\| \\
&\leq & \|T^*(I-\pi_n)\|\, \|(I-\pi_n)T\|\\
& =&  \|(I-\pi_n)T\|^2.
\eeqarray
Thus, we obtain a better estimate for  $\|T^*T- T_n^*T_n\|$ than the one given in Theorem \ref{Th-special-1}. 
}\erem

\bt\label{Th-special-2}
Suppose 
\ben
\i[(i)] $Z$ is a Banach space with a norm $\|\cdot\|_Z$ stonger than the norm on $Y$,  
\i[(ii)] $T: X\to Z$ is  a compact operator, 
\i[(iii)]  $Y_n\subseteq Z$ for every $n\in \N$,
\i[(iv)] $\pi_n: Z\to Z$ is a  bounded linear operator  such that 
$\|u-\pi_nu\|_Z\to 0$  {as}   $n\to\infty$
for every $u\in Z$.
\een
Then, $\|(I -\pi_n)T\|   \to 0$ {as} $n\to\infty.$  In particular,    $\|T^*T- T_n^*T_n\|\to 0$ {as} $n\to\infty.$
\et

\bpf
Under the given assumptions on $T$ and $(\pi_n)$, it follows  from Proposition \ref{Prop-fa}  that 
$$\|(I-\pi_n)T\|_{X\to Z}\to 0\q\h{as}\q n\to\infty,$$ 
where $\|\cdot\|_{X\to Z}$ is the norm on the space of bounded lonear operators from $X$ to $Z$.  Since the Banach space norm $\|\cdot\|_Z$ on $Z$ is stronger than the norm $\|\cdot\|_Y$ on $Y$, we obtain 
$$\|(I-\pi_n)T\| \leq c\,\|(I-\pi_n)T\|_{X\to Z}$$
for some constant $c>0$, so that $\|(I-\pi_n)T\| \to 0$ as $n\to\infty$. Hence, arrive at the conclusion by using Theorem \ref{Th-special-1}.
\epf

\section{Estimate under noisy data}

Suppose the availbale data is noisy, say we have $\tilde y\in Y$ such that 
$$\|\pi_n(y-\tilde y)\|\leq \d_n$$
for some $\d_n>0$.
Let 
$$\tilde x_n^\dagger = T_n^\dagger \pi_n\tilde y.$$

\bt\label{Th-3}
Let $\l_n$ be the smallest nonzero eighenvalue of $T_nT_n^*$ and let $\s_n:=\sqrt{\l_n}$. Then 
$$\|x_n^\dagger - \tilde x_n^\dagger \| \leq \frac{\d_n}{\s_n}.$$
If $x^\dagger = \f(T^*T)u$  for some $u\in X$, where $\f: (0,\infty)\to (0, \infty)$ is a continuous concave function, and if 
$$\d_n\leq \s_n\f(\e_n),$$
then 
$$\|x^\dagger - \tilde x_n^\dagger \| \leq c\f(\e_n),
$$
for some constant $c$, depending on $x^\dagger$.
\et

\bpf
Note that  
$$x_n^\dagger - \tilde x_n^\dagger = T_n^\dagger\pi_n (y-\tilde y).$$
Since nonzero eigenvalues of $T_n^*T_n$ and $T_nT_n^*$ are the same,  $\s_n:=\sqrt{\l_n}$ is the smallest singular value of $T_n$, so that  we have  (see (\cite{mtn-fa})
 $$\|T_n^\dagger\| = \frac{1}{\s_n}.$$
Thus, 
$$\|x_n^\dagger - \tilde x_n^\dagger \| \leq \frac{\d_n}{\s_n},$$
Using the assumption $x^\dagger = \f(T^*T)u$\,  for some $u\in X$, where $\f: (0,\infty)\to (0, \infty)$ is a continuous concave function, we know that 
$$\|x^\dagger - x_{\e_n}\| \leq \|u\|\f(\e_n).$$
Now, the last estimate in the theorem follows using the estimate 
$\|x^\dagger -  x_n^\dagger \|\leq 2 \|x^\dagger-x_{\e_n}\|$ obtained in Theorem \ref{Th-1}. 
\epf

\section{Further Regularization}

Recall that, in order to find $x_n^\dagger$ it is enough to solve the equation  (\ref{eq-comp-1}), that is, 
$$ T_nT_n^* v_n =  y_n $$
and take $x_n^\dagger = T_n^*v_n.$
However, under the noisy data  $\tilde y_n$ in place of $y_n:=P_ny$, it is not advisable to find $\tilde x_n^\dagger := T_n^\dagger \tilde y_n$; instead, we may solve a regularized equation. 

%

\bt\label{Th-5}
Let $ v_{\a,n}$ and   $\tilde v_{\a,n}$ be in $Y_n$ be such that 
$$(T_nT_n^* + \a I_n)v_{\a,n} =  y_n,\q (T_nT_n^* + \a I_n)\tilde v_{\a,n} = \tilde y_n$$ 
and let 
$$x_{\a,n}:= T_n^*v_{\a,n},\q \tilde x_{\a,n}:= T_n^* \tilde v_{\a,n}.$$
Then
$$\|x^\dagger - x_{\a,n}\| \leq \Big(1+\frac{\e_n}{\a}\Big)\|x^\dagger-x_\a\| $$
and 
$$\|x^\dagger - \tilde x_{\a,n}\| \leq \Big(1+\frac{\e_n}{\a}\Big)\|x^\dagger-x_\a\| + \frac{\|y_n-\tilde y_n\|}{\sqrt \a}. $$
In particular,  
$$\|x^\dagger - x_{\e_n,n}\| \leq 2 \|x^\dagger-x_{\e_n}\| $$
and 
$$\|x^\dagger - \tilde x_{\e_n,n}\| \leq 2 \|x^\dagger-x_{\e_n}\| + \frac{\|y_n-\tilde y_n\|}{\sqrt \e_n}. $$
\et

\bpf
From the definiton of  $v_{\a,n}$ and $\tilde v_{\a,n}$, it follows that 
$$ x_{\a,n}:= T_n^*  v_{\a,n},\q \tilde x_{\a,n}:= T_n^* \tilde v_{\a,n}$$
satisfy 
$$(T_n^*T_n + \a I)x_{\a,n}=T_n^* y_n,\q (T_n^*T_n + \a I)\tilde x_{\a,n}=T_n^* \tilde y_n.$$ 
Therefore,
$$\| x_{\a,n}- \tilde x_{\a,n}\|\leq \frac{\|y_n-\tilde y_n\|}{\sqrt\a} $$
Note also that 
\beqarray
x_\a-x_{\a, n} &=& (T^*T + \a I)^{-1} T^*y - (T_n^*T_n + \a I)^{-1} T_n^*y_n\\ 
&=& (T^*T + \a I)^{-1} T^*Tx^\dagger  - (T_n^*T_n + \a I)^{-1} T_n^*T_n x^\dagger\\
&=&T^*T (T^*T + \a I)^{-1} x^\dagger  - (T_n^*T_n + \a I)^{-1} T_n^*T_n x^\dagger\\
&=&    \a  (T_n^*T_n + \a I)^{-1} [ T^*T  -   T_n^*T_n  ] (T^*T + \a I)^{-1} x^\dagger
\eeqarray
Hence,
$$\| x_\a-x_{\a, n}\| \leq \e_n \|(T^*T + \a I)^{-1} x^\dagger\|.$$
But,
$$x^\dagger - x_{\a }= \a  (T^*T + \a I)^{-1}x^\dagger.$$
Therefore,
$$\|x_{\a,n} - x_\a\|\leq \frac{\e_n}{\a}\|x^\dagger-x_\a\|$$
so that 
$$\|x^\dagger-x_{\a, n}\|\leq \Big(1+\frac{\e_n}{\a}\Big)\|x^\dagger-x_\a\|.$$
Thus, we obtain the results. 
\epf

The following theorem is a consequence of Theorems \ref{Th-Tikh} and \ref{Th-5}.

\bt
Suppose $x^\dagger$ belongs to the range of $\f(T^*T)$, where  $\f: (0, \infty)\to (0, \infty)$ is a continuous  function  as in Theorem \ref{Th-Tikh}.
If  $\|y_n-\tilde y_n\|\leq \sqrt{\e_n}\f(\e_n),$ then 
$$\|x^\dagger -\tilde  x_{\e_n,n} \| \leq c\,\f(\e_n).$$
for some constant $c>0$ depending on $x^\dagger$.\et

\section{Applications to Integral Equations of the First Kind}

Let $\O:=[a, b]$, and $T$ be the Fredholm  integral operator on the real Hilbert space $L^2(\O)$ with a continuous kernel, that is,
$$(Tx)(s) = \int_a^b k(s, t) x(t)\, dt,\q x\in L^2(\O),\, s\in \O,$$
where $k\in C(\O\times \O)$.  Let $Z=C(\O)$ be with $\|\cdot\|_\infty$.

\beg\label{eg-interpol} \rm 
Let  $(\pi_n)$ be a sequence of  interpolatory projections with $R(\pi_n)\subseteq C(\O)$ and such that $(\pi_n)$ converges pointwise to the identity operator on $C(\O)$, that is, for each $u\in C(\O)$,
$$\|u-\pi_nu\|_\infty\to 0\q\h{as}\q n\to\infty.$$
In this case, $Y_n=R(\pi_n)$. 
Since $R(T)\subseteq C(\O)$  and $T: L^2[a, b]\to C(\O)$   is also a compact operator,  Proposition \ref{Prop-fa} and Theorem \ref{Th-special-2}  can be applied.

Note that  the  intepolatory projection  $\pi_n: C(\O)\to C(\O)$ can be representated as
$$
\pi_nx = \sum_{i=1}^n x(t_i^{(n)}) u_i^{(n)},\q x\in C(\O),$$
where  $t_1^{(n)}, \ldots, t_n^{(n)}$ are in $\O$ with  $t_1^{(n)}< \ldots < t_n^{(n)}$  and  $u_1^{(n)}, \ldots, u_n^{(n)}$ are in $C(\O)$ are such that 
$$u_i^{(n)}(t_j^{(n)}) = \left\{\begin{array}{ll} 1 & \h{if}\, i=j,\\  
0 & \h{if}\, i\not=j, \end{array}\right.,\q i, j=1, \ldots, n.$$
In this case,  for $x, v\in L^2(\O)$, we have 
\beqarray
\<T_nx, v\> &=& \<\pi_nTx, v\>\\
&=& \int_\O \sum_{i=1}^n (Tx)(t_i^{(n)})u_i^{(n)}(\t)  v(\t) d\t\\
&=& \int_\O \sum_{i=1}^n\(\int_\O k(t_i^{(n)}, t)x(t) dt  \)u_i^{(n)}(\t)  v(\t) d\t\\
&=& \int_\O x(t)\( \sum_{i=1}^n  k(t_i^{(n)}, t)\int_\O  u_i^{(n)}(\t)  v(\t) d\t\) dt.
\eeqarray
Thus,
$$(T_n^*v)(t) = \sum_{i=1}^n  k(t_i^{(n)}, t)\<u_i^{(n)},  v\>.$$
Hence,
\beqarray
T_nT_n^*v &=& \pi_nT(T_n^*v)\\
&=& \sum_{i=1}^n  (TT_n^*v)(t_i^{(n)}) u_i^{(n)}\\
&=& \sum_{i=1}^n \( \int_\O k(t_i^{(n)}, t) (T_n^*v)(t) dt \)u_i^{(n)}\\
&=& \sum_{i=1}^n \( \int_\O k(t_i^{(n)}, t)  \sum_{\ell=1}^n  k(t_\ell^{(n)}, t)\<u_\ell^{(n)},  v\> dt \)u_i^{(n)}
\eeqarray
Let  $\tilde y_n=\sum_{i=1}^n \b_i u_i$ and 
$$a_{ij} = \int_\O k(t_i^{(n)}, t) \(\sum_{\ell=1}^n k(t_{\ell}^{(n)}, t)\int_\O  u_\ell^{(n)}(\t) u_j^{(n)}(\t)d\t \) dt.$$
Then it can be seen that  $\tilde v_{\a,n}:= \sum_{i=1}^n x_i u_i$ is the solution of 
$$(T_nT_n^* + \a I_n)\tilde v_{\a,n} = \tilde y_n$$ 
if and only if   ${\bf x} = [x_i]$ is the solution of$${\bf A  x + \a x=  b},$$
  with  ${\bf A} = [a_{ij}]$ and ${\bf b}=[\b_i].$
\eeg 

\beg{\rm 
Let   $t_1^{(n)}, \ldots, t_n^{(n)}$ in $\O$ and $w_1^{(n)}, \ldots, w_n^{(n)}$ be positive real numbers such that the associated quadrature formula converges, that is, for every $f\in C(\O)$, 
$$\sum_{i=1}^n f(t_i^{(n)})w_i^{(n)}\to \int_a^b f(t)\, dt\q\h{as}\q n\to\infty.$$
Let  $Y_n:=\R^n_w$ be with inner product 
$$\<u, v\>_w:=\sum_{i=1}^n w_i^{(n)}u_i\bar v_i$$
for $u:=(u_i),\, v:=(v_i)$ in $\K^n.$
Let $\pi_n: C(\O)\to \K^n_w$ be defined by  
$$\pi_n f = (f(t_1^{(n)}), \ldots, f(t_1^{(n)})),\q f\in C(\O).$$
Thus, denoting the $i$-th coordinate of $u\in \R^n$ by $u(i)$,
$$(T_nx)(i) = (Tx)(t_i^{(n)}), \q i\in \{1, \ldots, n\}.$$
It is proved in   \cite{mtn-acm}   that the adjoint  
$T_n^*: \R_w^n\to L^2(\O)$ of $T_n$ is given by
$$(T_n^* \vec{\a})(s) = \sum_{i=1}^n k(\t_i^{(n)}, s)\a_i w_i^{(n)}, \q \vec{\a} = (\a_1, \ldots, \a_n)\in \R_w^n $$
so that 
$T_n^*T_n=F_nT$, where  $F_n$ is  the Nystr\"{o}m approximation of the integral operator $T^*: C(\O)\to C(\O)$ given by
$$(T^*u)(s) = \int_a^b {k(t, s)} u(t)\, dt,\q u\in C(\O),\, s\in \O,$$
that is,
$$(F_nx)(s) = \sum_{i=1}^n {k(t_i^{(n)}, s)} x(t_i^{(n)})w_i^{(n)},$$
%
%
and consequently, using the fact that  (see, e.g., \cite{kress, mtn-linop}) 
$\|F_nu-T^*u\|_\infty\to 0$  {as}  $n\to\infty$ and Proposition \ref{Prop-fa},  we have
$$\|T^*T-T_n^*T_n\| = \|(F_n-T^*)T\|\to 0\q\h{as}\q n\to\infty.$$ 
Thus, the method in \cite{mtn-acm}  is a special case.
In this case,  for $\vec{\a} = (\a_1, \ldots, \a_n)\in \R_w^n,$\, 
\beqarray
[(T_nT_n^*) \vec{\a}](i)  
&=& [T (T_n^* \vec{\a})])(t_i^{(n)})\\
&=& \int_\O k((t_i^{(n)}, t)[ (T_n^* \vec{\a})])(t) ] dt\\
&=& \sum_{j=1}^n \( w_j^{(n)}\int_\O k((t_i^{(n)}, t) k(\t_j^{(n)}, t)  dt\)\a_j 
\eeqarray
Then, writing 
$$a_{ij}:= w_j^{(n)}\int_\O k((t_i^{(n)}, t) k(\t_j^{(n)}, t)  dt\q\h{and}\q \tilde y_n:= (\b_1, \ldots, \b_n),$$ 
it can be seen that  $\tilde v_{\a,n}:= (x_1, \ldots, x_n)\in \R^n$  is the solution of 
$$(T_nT_n^* + \a I_n)\tilde v_{\a,n} = \tilde y_n$$ 
if and only if   ${\bf x} = [x_i]$ is the solution of$${\bf A  x + \a x=  b},$$
with  ${\bf A} = [a_{ij}]$ and $ {\bf b}=[\b_i].$ 
}\eeg

\vsq
\noi{\bf Acknowledgement.}  This work was done while I  was a Visiting Professor at Institute Camille Jordan, Saint Etienne for a month in June 2016.  I gratefully acknowledge the support  from the Institute Camille Jordan and the warm hospitality received from Prof. Mario Ahues and Prof. Laurence Grammont.

\end{document}